\font\fifteenrm=cmr10 scaled\magstep2 
\font\fifteeni=cmmi10 scaled\magstep2
\font\fifteensy=cmsy10 scaled\magstep2
\font\fifteenbf=cmbx10 scaled\magstep2
\font\fifteentt=cmtt10 scaled\magstep2
\font\fifteenit=cmti10 scaled\magstep2
\font\fifteensl=cmsl10 scaled\magstep2
\font\fifteenam=msam10 scaled\magstep2
\font\fifteenbm=msbm10 scaled\magstep2
\font\fifteenex=cmex10 scaled\magstep2
\font\fifteensc=cmcsc10 scaled\magstep2 
\font\twelverm=cmr10 at 12pt
\font\twelvei=cmmi10 at 12pt
\font\twelvesy=cmsy10 at 12pt
\font\twelvebf=cmbx10 at 12pt
\font\twelvett=cmtt10 at 12pt
\font\twelveit=cmti10 at 12pt
\font\twelvesl=cmsl10 at 12pt
\font\twelveam=msam10 at 12pt
\font\twelvebm=msbm10 at 12pt
\font\twelveex=cmex10 at 12pt
\font\twelvesc=cmcsc10 at 12pt
\font\elevenrm=cmr10 scaled\magstephalf 
\font\eleveni=cmmi10 scaled\magstephalf
\font\elevensy=cmsy10 scaled\magstephalf
\font\elevenbf=cmbx10 scaled\magstephalf
\font\eleventt=cmtt10 scaled\magstephalf
\font\elevenit=cmti10 scaled\magstephalf
\font\elevensl=cmsl10 scaled\magstephalf
\font\elevenam=msam10 scaled\magstephalf
\font\elevenbm=msbm10 scaled\magstephalf
\font\elevenex=cmex10 scaled\magstephalf
\font\elevensc=cmcsc10 scaled\magstephalf
\font\tenrm=cmr10
\font\teni=cmmi10
\font\tensy=cmsy10
\font\tenbf=cmbx10
\font\tentt=cmtt10
\font\tenit=cmti10
\font\tensl=cmsl10
\font\tenam=msam10
\font\tenbm=msbm10
\font\tenex=cmex10
\font\tensc=cmcsc10
\font\ninerm=cmr9
\font\ninei=cmmi9
\font\ninesy=cmsy9
\font\ninebf=cmbx9
\font\ninett=cmtt9
\font\nineit=cmti9
\font\ninesl=cmsl9
\font\nineam=msam9
\font\ninebm=msbm9
\font\nineex=cmex9
\font\ninesc=cmcsc9
\font\eightrm=cmr8
\font\eighti=cmmi8
\font\eightsy=cmsy8
\font\eightbf=cmbx8
\font\eighttt=cmtt8
\font\eightit=cmti8
\font\eightsl=cmsl8
\font\eightam=msam8
\font\eightbm=msbm8
\font\eightex=cmex8
\font\eightsc=cmcsc8
\font\sevenrm=cmr7
\font\seveni=cmmi7
\font\sevensy=cmsy7
\font\sevenbf=cmbx7

\font\sevenam=msam7
\font\sevenbm=msbm7

\font\sixrm=cmr6
\font\sixi=cmmi6
\font\sixsy=cmsy6

\font\sixam=msam6
\font\sixbm=msbm6

\font\fiverm=cmr5
\font\fivei=cmmi5
\font\fivesy=cmsy5

\font\fiveam=msam5
\font\fivebm=msbm5

\font\fourrm=cmr5 at 4pt
\font\fouri=cmmi5 at 4pt
\font\foursy=cmsy5 at 4pt

\font\fouram=msam5 at 4pt
\font\fourbm=msbm5 at 4pt

\skewchar\twelvei='177 \skewchar\eleveni='177\skewchar\teni='177
\skewchar\ninei='177 \skewchar\eighti='177\skewchar\seveni='177 
\skewchar\sixi='177 \skewchar\fivei='177 \skewchar\fouri='177
\skewchar\twelvesy='60 \skewchar\elevensy='60 \skewchar\tensy='60
\skewchar\ninesy='60 \skewchar\eightsy='60 \skewchar\sevensy='60 
\skewchar\sixsy='60 \skewchar\fivesy='60 \skewchar\foursy='60
\newfam\itfam
\newfam\slfam
\newfam\bffam
\newfam\ttfam
\newfam\scfam
\newfam\amfam
\newfam\bmfam
\def\eightbig#1{{\hbox{$\left#1\vbox to 6.5pt{}\voidright $}}}
\def\eightBig#1{{\hbox{$\left#1\vbox to 7.5pt{}\voidright $}}}
\def\eightbigg#1{{\hbox{$\left#1\vbox to 10pt{}\voidright $}}}
\def\eightBigg#1{{\hbox{$\left#1\vbox to 13pt{}\voidright $}}}
\def\ninebig#1{{\hbox{$\left#1\vbox to 7.5pt{}\voidright $}}}
\def\nineBig#1{{\hbox{$\left#1\vbox to 8.5pt{}\voidright $}}}
\def\ninebigg#1{{\hbox{$\left#1\vbox to 11.5pt{}\voidright $}}}
\def\nineBigg#1{{\hbox{$\left#1\vbox to 14.5pt{}\voidright $}}}
\def\tenbig#1{{\hbox{$\left#1\vbox to 8.5pt{}\voidright $}}}
\def\tenBig#1{{\hbox{$\left#1\vbox to 9.5pt{}\voidright $}}}
\def\tenbigg#1{{\hbox{$\left#1\vbox to 12.5pt{}\voidright $}}}
\def\tenBigg#1{{\hbox{$\left#1\vbox to 16pt{}\voidright $}}}
\def\elevenbig#1{{\hbox{$\left#1\vbox to 9pt{}\voidright $}}}
\def\elevenBig#1{{\hbox{$\left#1\vbox to 10.5pt{}\voidright $}}}
\def\elevenbigg#1{{\hbox{$\left#1\vbox to 14pt{}\voidright $}}}
\def\elevenBigg#1{{\hbox{$\left#1\vbox to 17.5pt{}\voidright $}}}
\def\twelvebig#1{{\hbox{$\left#1\vbox to 10pt{}\voidright $}}}
\def\twelveBig#1{{\hbox{$\left#1\vbox to 11pt{}\voidright $}}}
\def\twelvebigg#1{{\hbox{$\left#1\vbox to 15pt{}\voidright $}}}
\def\twelveBigg#1{{\hbox{$\left#1\vbox to 19pt{}\voidright $}}}
\def\fifteenbig#1{{\hbox{$\left#1\vbox to 12pt{}\voidright $}}}
\def\fifteenBig#1{{\hbox{$\left#1\vbox to 13.5pt{}\voidright $}}}
\def\fifteenbigg#1{{\hbox{$\left#1\vbox to 18pt{}\voidright $}}}
\def\fifteenBigg#1{{\hbox{$\left#1\vbox to 23pt{}\voidright $}}}
\def\voidright{\right.\nulldelimiterspace=0pt \mathsurround=0pt }
\def\fifteenpoint{
  \textfont0=\fifteenrm \scriptfont0=\twelverm \scriptscriptfont0=\tenrm
  \def\rm{\fam0 \fifteenrm}%
  \textfont1=\fifteeni \scriptfont1=\twelvei \scriptscriptfont1=\teni
  \textfont2=\fifteensy \scriptfont2=\twelvesy \scriptscriptfont2=\tensy
  \textfont3=\fifteenex \scriptfont3=\fifteenex \scriptscriptfont3=\fifteenex
  \def\it{\fam\itfam\fifteenit}\textfont\itfam=\fifteenit
  \def\sl{\fam\slfam\fifteensl}\textfont\slfam=\fifteensl
  \def\bf{\fam\bffam\fifteenbf}\textfont\bffam=\fifteenbf 
    \scriptfont\bffam=\twelvebf\scriptscriptfont\bffam=\tenbf
  \def\tt{\fam\ttfam\fifteentt}\textfont\ttfam=\fifteentt
  \def\sc{\fam\scfam\fifteensc}\textfont\scfam=\fifteensc
  \def\am{\fam\amfam\fifteenam}\textfont\amfam=\fifteenam
    \scriptfont\amfam=\twelveam\scriptscriptfont\amfam=\tenam
  \def\bm{\fam\bmfam\fifteenbm}\textfont\bmfam=\fifteenbm
    \scriptfont\bmfam=\twelvebm\scriptscriptfont\bmfam=\tenbm
  \baselineskip=21pt \rm
  \let\big=\fifteenbig\let\Big=\fifteenBig\let\bigg=\fifteenbigg
  \let\Bigg=\fifteenBigg}
\def\twelvepoint{
  \textfont0=\twelverm \scriptfont0=\ninerm \scriptscriptfont0=\sevenrm
  \def\rm{\fam0 \twelverm}%
  \textfont1=\twelvei \scriptfont1=\ninei \scriptscriptfont1=\seveni
  \textfont2=\twelvesy \scriptfont2=\ninesy \scriptscriptfont2=\sevensy
  \textfont3=\twelveex \scriptfont3=\twelveex \scriptscriptfont3=\twelveex
  \def\it{\fam\itfam\twelveit}\textfont\itfam=\twelveit
  \def\sl{\fam\slfam\twelvesl}\textfont\slfam=\twelvesl
  \def\bf{\fam\bffam\twelvebf}\textfont\bffam=\twelvebf 
    \scriptfont\bffam=\ninebf\scriptscriptfont\bffam=\sevenbf
  \def\tt{\fam\ttfam\twelvett}\textfont\ttfam=\twelvett
  \def\sc{\fam\scfam\twelvesc}\textfont\scfam=\twelvesc
  \def\am{\fam\amfam\twelveam}\textfont\amfam=\twelveam
    \scriptfont\amfam=\nineam\scriptscriptfont\amfam=\sevenam
  \def\bm{\fam\bmfam\twelvebm}\textfont\bmfam=\twelvebm
    \scriptfont\bmfam=\ninebm\scriptscriptfont\bmfam=\sevenbm
  \baselineskip=17.8pt \rm 
  \def\looselineskip{\baselineskip=18.5pt plus 1.8pt}%
  \def\tightlineskip{\baselineskip=16.5pt}%
  \def\verytightlineskip{\baselineskip=15pt}%
  \let\big=\twelvebig\let\Big=\twelveBig\let\bigg=\twelvebigg
  \let\Bigg=\twelveBigg  }
\def\elevenpoint{
  \textfont0=\elevenrm \scriptfont0=\ninerm \scriptscriptfont0=\sixrm
  \def\rm{\fam0 \elevenrm}%
  \textfont1=\eleveni \scriptfont1=\ninei \scriptscriptfont1=\sixi
  \textfont2=\elevensy \scriptfont2=\ninesy \scriptfont2=\sixsy 
  \textfont3=\elevenex \scriptfont3=\elevenex \scriptfont3=\elevenex
  \def\it{\fam\itfam\elevenit}\textfont\itfam=\elevenit
  \def\sl{\fam\slfam\elevensl}\textfont\slfam=\elevensl
  \def\bf{\fam\bffam\elevenbf}\textfont\bffam=\elevenbf
  \def\tt{\fam\ttfam\eleventt}\textfont\ttfam=\eleventt
  \def\sc{\fam\scfam\elevensc}\textfont\scfam=\elevensc
  \def\am{\fam\amfam\elevenam}\textfont\amfam=\elevenam
    \scriptfont\amfam=\nineam\scriptscriptfont\amfam=\sixam
  \def\bm{\fam\bmfam\elevenbm}\textfont\bmfam=\elevenbm
    \scriptfont\bmfam=\ninebm\scriptscriptfont\bmfam=\sixbm
  \baselineskip=15.1pt \rm
  \def\looselineskip{\baselineskip=16pt plus 1.5pt}%
  \def\tightlineskip{\baselineskip=14pt}%
  \def\verytightlineskip{\baselineskip=13pt}%
  \let\big=\elevenbig\let\Big=\elevenBig\let\bigg=\elevenbigg
  \let\Bigg=\elevenBigg  }
\def\tenpoint{
  \textfont0=\tenrm \scriptfont0=\eightrm \scriptscriptfont0=\fiverm
  \def\rm{\fam0 \tenrm}%
  \textfont1=\teni \scriptfont1=\eighti \scriptscriptfont1=\fivei
  \textfont2=\tensy \scriptfont2=\eightsy \scriptfont2=\fivesy 
  \textfont3=\tenex \scriptfont3=\tenex \scriptfont3=\tenex
  \def\it{\fam\itfam\tenit}\textfont\itfam=\tenit
  \def\sl{\fam\slfam\tensl}\textfont\slfam=\tensl
  \def\bf{\fam\bffam\tenbf}\textfont\bffam=\tenbf
  \def\tt{\fam\ttfam\tentt}\textfont\ttfam=\tentt
  \def\sc{\fam\scfam\tensc}\textfont\scfam=\tensc
  \def\am{\fam\amfam\tenam}\textfont\amfam=\tenam
    \scriptfont\amfam=\eightam \scriptscriptfont\amfam=\fiveam
  \def\bm{\fam\bmfam\tenbm}\textfont\bmfam=\tenbm
    \scriptfont\bmfam=\eightbm \scriptscriptfont\bmfam=\fivebm
  \baselineskip=14pt\rm
  \def\looselineskip{\baselineskip=14.8pt plus1.5pt}
  \def\tightlineskip{\baselineskip=12.6pt}%
  \def\verytightlineskip{\baselineskip=13pt}%
  \let\big=\tenbig\let\Big=\tenBig\let\bigg=\tenbigg\let\Bigg=\tenBigg  }
\def\ninepoint{
  \textfont0=\ninerm \scriptfont0=\sevenrm \scriptscriptfont0=\fourrm
  \def\rm{\fam0 \ninerm}%
  \textfont1=\ninei \scriptfont1=\seveni \scriptscriptfont1=\fouri
  \textfont2=\ninesy \scriptfont2=\sevensy \scriptfont2=\foursy 
  \textfont3=\nineex \scriptfont3=\nineex \scriptfont3=\nineex
  \def\it{\fam\itfam\nineit}\textfont\itfam=\nineit
  \def\sl{\fam\slfam\ninesl}\textfont\slfam=\ninesl
  \def\bf{\fam\bffam\ninebf}\textfont\bffam=\ninebf
  \def\tt{\fam\ttfam\ninett}\textfont\ttfam=\ninett
  \def\sc{\fam\scfam\ninesc}\textfont\scfam=\ninesc
  \def\am{\fam\amfam\nineam}\textfont\amfam=\nineam
    \scriptfont\amfam=\nineam\scriptscriptfont\amfam=\fouram
  \def\bm{\fam\bmfam\ninebm}\textfont\bmfam=\ninebm
    \scriptfont\bmfam=\ninebm\scriptscriptfont\bmfam=\fourbm
  \baselineskip=12.6pt\rm
  \def\tightlineskip{\baselineskip=11.5pt}
  \let\big=\ninebig\let\Big=\nineBig\let\bigg=\ninebigg
  \let\Bigg=\nineBigg  }
\def\eightpoint{
  \textfont0=\eightrm \scriptfont0=\fiverm \scriptscriptfont0=\fourrm
  \def\rm{\fam0 \eightrm}%
  \textfont1=\eighti \scriptfont1=\fivei \scriptscriptfont1=\fouri
  \textfont2=\eightsy \scriptfont2=\fivesy \scriptfont2=\foursy 
  \textfont3=\eightex \scriptfont3=\eightex \scriptfont3=\eightex
  \def\it{\fam\itfam\eightit}\textfont\itfam=\eightit
  \def\sl{\fam\slfam\eightsl}\textfont\slfam=\eightsl
  \def\bf{\fam\bffam\eightbf}\textfont\bffam=\eightbf
  \def\tt{\fam\ttfam\eighttt}\textfont\ttfam=\eighttt
  \def\sc{\fam\scfam\eightsc}\textfont\scfam=\eightsc
  \def\am{\fam\amfam\eightam}\textfont\amfam=\eightam
    \scriptfont\amfam=\eightam\scriptscriptfont\amfam=\fouram
  \def\bm{\fam\bmfam\eightbm}\textfont\bmfam=\eightbm
    \scriptfont\bmfam=\eightbm\scriptscriptfont\bmfam=\fourbm
  \baselineskip=11.2pt \rm
  \def\tightlineskip{\baselineskip=10.4pt}
  \let\big=\eightbig\let\Big=\eightBig\let\bigg=\eightbigg
  \let\Bigg=\eightBigg  }

\twelvepoint
\nopagenumbers
\hsize=6in\vsize=8.8in

\parskip=1pt plus 1pt

\newif\ifSpecialhead\Specialheadfalse
\newbox\specialheadbox

\def\specialhead #1\par{\Specialheadtrue\setbox\specialheadbox=\hbox{#1}}
\headline={{\ifSpecialhead\box\specialheadbox\global\Specialheadfalse\else
     \ifnum\pageno<0{\hfill\quad{\twelvebf\folio}}%
     \else\ifnum\pageno<2\hfill
     \else\hfill\twelvepoint\sc\firstmark\quad{\twelvebf\folio}\fi\fi\fi}}

\def\title#1\par{\bigskip{\def\cr{\par\center}\center\fifteenbf #1\par}\medskip}
\def\subtitle#1\par{\centerline{\fifteenrm #1}\medskip}
\def\author#1\par{\medskip{\def\cr{\par\center\twelvesc}\fifteensc\center#1\par}}
\def\center#1\par{\hfil #1\hfil\par}
\def\abstract.#1\par{\message{Abstract.}%
                    \medskip{\narrower\narrower\tenpoint\tightlineskip
                        \noindent{\bf Abstract.}#1\par}\medskip\noindent}
\def\tinyabstract.#1\par{\message{Abstract.}%
                    \medskip{\narrower\narrower\eightpoint\tightlineskip
                        \noindent{\bf Abstract.}#1\par}\medskip\noindent}
\def\bigabstract.#1\par{\message{Abstract.}%
                         \medskip{\narrower\narrower\tightlineskip
                         \noindent{\bf Abstract. }#1\par}\medskip\noindent}
\def\acknowledgement#1\par{\footnote{}{#1}}
\def\sectionskip{\Goodbreak\vskip 25pt plus 15pt minus 5pt}
\def\secnumber{\ifquiet
               \else\ifNoSections
                    \else\sectionsymbol\the\secno\quad\fi\fi}
\def\section#1\par{ \NoSectionsfalse\par\sectionskip\proofdepth=0\claimno=0
 \ifquiet\else\advance\secno by1\fi\toks0={#1}
 \immediate\write16{\ifquiet\else Section \the\secno\space\fi
                    \the\toks0}%
 \mark{\secnumber #1}%
 {\fifteenpoint\bf\noindent\secnumber #1}\nobreak\bigskip\quietoff
 \nobreak\noindent}
\def\quiet{\quiettrue}

\def\quietoff{\ifQUIET\else\quietfalse\fi}
\newif\ifquiet
\newif\ifQUIET
\newif\ifNoSections
\newcount\claimtype
\newcount\secno
\newcount\claimno
\newcount\subclaimno
\newcount\subsubclaimno
\newcount\subsubsubclaimno
\newcount\proofdepth
\def\subclaimnumber{\ifquiet\else\ifcase\subclaimno\or A\or B\or C\or D\or E\or
     F\or G\or H\or I\or J\or K\or L\or M\or N\or O\or P\fi\fi}
\def\subsubclaimnumber{\ifquiet\else\ifcase\subsubclaimno\or i\or ii\or iii\or 
   iv\or v\or vi\or vii\or viii\or ix\or x\or xi\or xii\or xiii\or xiv\fi\fi}
\def\subsubsubclaimnumber{\ifquiet\else\ifcase\subsubsubclaimno\or a\or b\or 
   c\or d\or e\or f\or g\or \or h\or i\or j\or k\or l\or m\or n\or o\fi\fi}
\def\claimtag{\ifquiet\else
  \ifNoSections
    \ifcase\proofdepth\the\claimno%
    \or\the\claimno.\subclaimnumber
    \or\the\claimno.\subclaimnumber.\subsubclaimnumber
    \or\the\claimno.\subclaimnumber.\subsubclaimnumber
                                                .\subsubsubclaimnumber\fi
  \else
    \ifcase\proofdepth\the\secno.\the\claimno
    \or\the\secno.\the\claimno.\subclaimnumber
    \or\the\secno.\the\claimno.\subclaimnumber.\subsubclaimnumber
    \or\the\secno.\the\claimno.\subclaimnumber.\subsubclaimnumber
                                                .\subsubsubclaimnumber\fi\fi\fi}
\secno=0\claimno=0\proofdepth=0\subclaimno=0\subsubclaimno=0\subsubsubclaimno=0
\NoSectionstrue
\newbox\qedbox
\def\claimname{\ifcase\claimtype Theorem\or Lemma\or Claim\or Corollary\or
               Question\or Definition\or Remark\or Conjecture\fi}
\def\preclaimskip{\removelastskip
    \ifcase\claimtype\goodbreak\vskip 8pt plus 4pt minus 2pt
                  \or\goodbreak\vskip 6pt plus 4pt minus 1pt
                  \or\goodbreak\vskip 5pt plus 4pt minus 1pt
                  \or\goodbreak\vskip 8pt plus 4pt minus 2pt
                  \or\vskip 7pt plus 4pt minus 2pt
                  \or\vskip 7pt plus 4pt minus 2pt
                  \or\vskip 7pt plus 4pt minus 2pt
                  \or\goodbreak\vskip 8pt plus 4pt minus 2pt\fi}
\def\postclaimskip{\ifcase\claimtype         \vskip 4pt plus 2pt minus 2pt
                                          \or\vskip 3pt plus 2pt minus 2pt
                                          \or\vskip 2pt plus 2pt minus 1pt
                                          \or\vskip 4pt plus 2pt minus 2pt
                                          \or\vskip 1pt plus 2pt 
                                          \or\vskip 4pt plus 4pt 
                                          \or\vskip 3pt plus 2pt
                                          \or\vskip 4pt plus 2pt minus 2pt\fi}
\def\claimfont{\ifcase\claimtype
                  \sl\or\sl\or\sl\or\sl\or\sl\or\rm\or\rm\or\sl\fi}
\def\advancetag{\ifcase\proofdepth\advance\claimno by1
                               \or\advance\subclaimno by1
                               \or\advance\subsubclaimno by1
                               \or\advance\subsubsubclaimno by1\fi}
\def\sayclaim#1.#2 #3\par{\ifquiet\else\advancetag\fi
    \preclaimskip\setbox1=\hbox{#1}\setbox2=\hbox{#2}%
    \toks0={#1 }
    \immediate\write16{\ifdim\wd1>0pt\the\toks0
                       \else\claimname\space\fi \claimtag.}%
    \vbox{\noindent
    {\bf\ifdim\wd1=0pt \claimname\else #1\fi\ifquiet.\else\ \claimtag{\ifNoSections.\fi}\fi}%
    \enspace{\ifdim\wd2>0pt\sc #2\enspace\fi}%
    {\claimfont #3\par}}\postclaimskip\quietoff}
\def\theorem{\claimtype=0\sayclaim}

\def\corollary{\claimtype=3\sayclaim}

\def\point#1. #2\par{\item{\rm #1.}#2\par}
\def\points#1\cr\par{\medskip\vbox{\let\cr=\point\point#1\par}\par}
\def\df{\it}
\def\prooffont{}
\def\proofsize{}
\def\proofindent{}
\def\proofskip{\badbreak\ifcase\claimtype    \vskip 3pt plus 2pt minus 2pt
                                          \or\vskip 2pt plus 2pt minus 2pt
                                          \or\vskip 1pt plus 2pt minus 1pt
                                          \or\vskip 3pt plus 2pt minus 2pt
                                          \or\vskip 1pt plus 2pt 
                                          \or\vskip 2pt plus 4pt 
                                          \or\vskip 1pt plus 2pt
                                          \or\vskip 3pt plus 2pt minus 2pt\fi}

\def\Goodbreak{\vskip0pt plus.5in\penalty-1000\vskip0pt plus-.5in}
\def\goodbreak{\penalty-500}
\def\badbreak{\penalty500}
\def\Badbreak{\penalty1000}
\def\proof{\message{proof}\removelastskip\Badbreak\proofskip\begingroup
  \advance\proofdepth by1
  \setbox\qedbox=\hbox{\halmos\raise2pt\hbox{\fiverm\claimname}}%
  \prooffont\proofsize\proofindent\noindent{\bf Proof: }}
\def\proofof#1:{\message{proof}\removelastskip\Badbreak\proofskip\begingroup
  \advance\proofdepth by1
  \setbox\qedbox=\hbox{\halmos\raise2pt\hbox{\fiverm#1}}%
  \prooffont\proofsize\proofindent\noindent{\bf Proof of #1: }}
\def\cite[#1]{[{\tenrm{#1}}]\message{[#1]}}
\edef\ref#1{\expandafter\global\expandafter\edef#1{\noexpand\claimtag}}
\newwrite\notes
\openout\notes=\jobname.notes
\long\def\unexpandedwrite#1#2{\def\finwrite{\write#1}%
   {\aftergroup\finwrite\aftergroup{\sanitize#2\endsanity}}}
\def\sanitize{\futurelet\next\sanswitch}
\let\stoken=\space
\def\sanswitch{\ifx\next\endsanity
  \else\ifcat\noexpand\next\stoken\aftergroup\space\let\next=\eat
   \else\ifcat\noexpand\next\bgroup\aftergroup{\let\next=\eat
    \else\ifcat\noexpand\next\egroup\aftergroup}\let\next=\eat
     \else\let\next=\copytoken\fi\fi\fi\fi \next}
\def\eat{\afterassignment\sanitize \let\next= }
\long\def\copytoken#1{\ifcat\noexpand#1\relax\aftergroup\noexpand
  \else\ifcat\noexpand#1\noexpand~\aftergroup\noexpand\fi\fi
  \aftergroup#1\sanitize}
\def\endsanity\endsanity{}

\def\note#1#2{\hbox to2in{\strut#1\quad\dotfill\quad#2}}
\def\boxit#1{\setbox4=\hbox{\kern1pt#1\kern1pt}
  \hbox{\vrule\vbox{\hrule\kern1pt\box4\kern1pt\hrule}\vrule}}
\def\halmos{\hbox{\am\char'3}} 
\def\qed#1\par{\message{.                                }\setbox1=\hbox{#1}%
  \ifdim\wd1>0pt\setbox\qedbox=\hbox{\halmos\raise2pt\hbox{\fiverm #1}}\fi
  \kern5pt\lower 2pt\hbox{\box\qedbox}\proofskip\goodbreak\endgroup}

\def\sectionsymbol{\S}
\def\k{\kappa}
\def\g{\gamma}
\def\a{\alpha}
\def\b{\beta}
\def\d{\delta}

\def\l{\lambda}
\def\z{\zeta}
\def\I1{\mathop{\hbox{\sc i}_1}}
\def\w{\omega}
\def\P{{\mathchoice{\hbox{\bm P}}{\hbox{\bm P}}
         {\hbox{\tenbm P}}{\hbox{\sevenbm P}}}}
\def\Q{{\mathchoice{\hbox{\bm Q}}{\hbox{\bm Q}}
         {\hbox{\tenbm Q}}{\hbox{\sevenbm Q}}}}

\def\X{{\mathchoice{\hbox{\bm X}}{\hbox{\bm X}}
         {\hbox{\tenbm X}}{\hbox{\sevenbm X}}}}
\def\card#1{\left|#1\right|}

\def\id{\mathop{\hbox{\tenrm id}}}

\def\elesub{\prec}

\def\unifto{\buildrel\lower 7pt\hbox{$\to$}\over\to}

\def\cof{\mathop{\rm cof}\nolimits}
\def\cp{\mathop{\rm cp}\nolimits}

\def\ORD{\hbox{\sc ord}}

\def\CH{\hbox{\sc ch}}

\def\in{\mathrel{\mathchoice{\raise 
1pt\hbox{$\scriptstyle\cal\char'62$}}
         {\raise 1pt\hbox{$\scriptstyle\cal\char'62$}}
         {\raise .5pt\hbox{$\scriptscriptstyle\cal\char'62$}}
         {\hbox{$\scriptscriptstyle\cal\char'62$}}}\penalty700{}}
\def\ni{\mathrel{\mathchoice{\raise 1pt\hbox{$\scriptstyle\cal\char'63$}}
                   {\raise 1pt\hbox{$\scriptstyle\cal\char'63$}}
                   {\raise .5pt\hbox{$\scriptscriptstyle\cal\char'63$}}
                   {\hbox{$\scriptscriptstyle\cal\char'63$}}}\penalty700}
\def\of{\mathrel{\mathchoice{\raise 1pt\hbox{$\scriptstyle\subseteq$}}
                   {\raise 1pt\hbox{$\scriptstyle\subseteq$}}
                   {\raise .5pt\hbox{$\scriptscriptstyle\subseteq$}}
                   {\hbox{$\scriptscriptstyle\subseteq$}}}}
\def\fo{\mathrel{\mathchoice{\raise 1pt\hbox{$\scriptstyle\supseteq$}}
                   {\raise 1pt\hbox{$\scriptstyle\supseteq$}}
                   {\raise .5pt\hbox{$\scriptscriptstyle\supseteq$}}
                   {\hbox{$\scriptscriptstyle\supseteq$}}}}
\def\notin{\mathrel{\mathchoice
  {\raise 1pt\hbox{\rlap{$\scriptstyle\;|$}$\scriptstyle\cal\char'62$}}
  {\raise 1pt\hbox{\rlap{$\scriptstyle\kern2pt 
          |$}$\scriptstyle\cal\char'62$}}
  {\raise .5pt\hbox{\rlap{$\scriptscriptstyle\, |$}$\scriptscriptstyle
      \cal\char'62$}}
  {\hbox{\rlap{$\scriptscriptstyle\, |$}$\scriptscriptstyle
     \cal\char'62$}}}%
  \penalty700}

\def\and{\mathrel{\kern1pt\&\kern1pt}}
\def\iff{\mathrel{\leftrightarrow}}

\def\union{\cup}

\def\intersect{\cap}

\def\cross{\times}

\def\[#1]{\left[\vphantom{\bigm|}#1\right]}
\def\<#1>{\langle\,#1\,\rangle}

\def\image{\mathbin{\hbox{\tt\char'42}}}
\def\restrict{\mathbin{\mathchoice{\hbox{\am\char'26}}{\hbox{\am\char'26}}{\hbox{\eightam\char'26}}{\hbox{\sixam\char'26}}}}
\def\force{\mathbin{\hbox{\am\char'15}}}

\def\st{\mid}
\def\seq<#1>{{\def\st{\mid\penalty650}\left<\,#1\,\right>}}

\def\set#1{\{\,#1\,\}}

\def\th{{\hbox{\fiverm th}}}

\def\forces{\force}
\def\lttheta{{\raise 1pt\hbox{$\scriptstyle<$}\theta}}

\def\I1{\mathop{\hbox{\sc i}_1}}
\def\ltk{{{\scriptstyle<}\k}}
\def\ltl{{{\scriptstyle<}\l}}

\def\lteb{{{\scriptstyle\leq}\b}}
\def\lted{{{\scriptstyle\leq}\d}}

\def\Qdot{\dot\Q}

\def\Qdot{\dot\Q}

\def\jmu{j_\mu}
\def\jnu{j_\nu}

\def\Vbar{{\overline V}}
\def\Mbar{{\overline M}}

\font\arrow=line10 scaled \magstep1
\def\makeline#1.{\hbox{\arrow\char#1}}
\def\makearrow#1.#2.{\hbox{\arrow\char#1\llap{\char#2}}}
\def\definelinesandarrows#1.#2.#3.#4.#5.{
   \expandafter\edef\csname#4line\endcsname{\makeline#1.}
   \expandafter\edef\csname#4arrow\endcsname{\makearrow#1.#2.}
   \expandafter\edef\csname#5line\endcsname{\makeline#1.}
   \expandafter\edef\csname#5arrow\endcsname{\makearrow#1.#3.}}
\definelinesandarrows 0.18.9.ne.sw.
\definelinesandarrows 1.21.11.nnne.sssw.
\definelinesandarrows 2.14.13.nnnne.ssssw.
\definelinesandarrows 3.23.15.nnnnne.sssssw.
\definelinesandarrows 4.23.15.nnnnnne.ssssssw.
\definelinesandarrows 10.30.29.nne.ssw.
\definelinesandarrows 16.49.41.neeeeee.swwwwww.
\definelinesandarrows 17.51.43.neeee.swwww.
\definelinesandarrows 19.55.47.nehuh.swhuh.
\definelinesandarrows 24.58.41.neeeeeee.swwwwwww.
\definelinesandarrows 26.62.9.neee.swww.
\definelinesandarrows 33.49.25.neeeee.swwwww.
\definelinesandarrows 35.62.61.nee.sww.
\definelinesandarrows 64.82.73.se.nw.
\definelinesandarrows 65.85.75.ssse.nnnw.
\definelinesandarrows 66.78.77.sssse.nnnnw.
\definelinesandarrows 67.87.79.ssssse.nnnnnw.
\definelinesandarrows 68.87.79.sssssse.nnnnnnw.
\definelinesandarrows 74.94.93.sse.nnw.
\definelinesandarrows 80.113.105.seeeeee.nwwwwww.
\definelinesandarrows 81.115.107.seeee.nwwww.
\definelinesandarrows 99.126.125.see.nww.
\def\sejoin#1#2{\setbox1=\hbox{#1}\setbox2=\hbox{#2}%
  \hbox{\vbox{\hbox{\copy1\kern\wd2}\nointerlineskip
              \hbox{\kern\wd1\box2}}}}
\def\nejoin#1#2{\setbox1=\hbox{#1}\setbox2=\hbox{#2}%
  \hbox{\vbox{\hbox{\kern\wd1\copy2}\nointerlineskip\hbox{\copy1\kern\wd2}}}}
\newdimen\hnudge
\newdimen\vnudge
\newdimen\hnudgedefault
\newdimen\vnudgedefault

\def\SEdefaultnudge{\hnudge=-16pt\vnudge=20pt}
\def\Edefaultnudge{\hnudge=-25pt\vnudge=6pt}

\def\longEdefaultnudge{\hnudge=-5pt\vnudge=6pt}
\def\nudgeright#1pt{\advance\hnudge by#1pt}
\def\nudgeleft#1pt{\advance\hnudge by-#1pt}
\def\nudgeup#1pt{\advance\vnudge by#1pt}
\def\nudgedown#1pt{\advance\vnudge by-#1pt}
\def\label#1{\smash{\llap{\kern\hnudge
                   \raise\vnudge\rlap{$\scriptstyle#1$}\hfill}}}

\def\SEarrow{\SEdefaultnudge
             \sejoin\seeline{\sejoin\seeline{\sejoin\seeline\seearrow}}}

\def\Earrow{\Edefaultnudge\setbox1=\hbox{\SEarrow}
 \hbox{\raise 2pt\hbox{\vrule height-.4pt depth.8ptwidth\wd1\kern2pt
       \llap{\arrow\char'55}}}}
\def\longEarrow{\longEdefaultnudge\setbox1=\hbox{\SEarrow}
      \rlap{\hskip-1.25\wd1\raise 2pt
            \hbox{\vrule height-.4pt depth.8ptwidth2.5\wd1\kern2pt
            \llap{\arrow\char'55}}}}

\looselineskip
\def\Mbar{{\overline M}}
\def\Vbar{{\overline V}}


\title Gap Forcing: generalizing the Levy-Solovay theorem

\author Joel David Hamkins\cr
	Kobe University and\cr
	The City University of New York\cr
	{\tentt http://www.math.csi.cuny.edu/$\sim$hamkins}\cr

\abstract. The Levy-Solovay Theorem \cite[LevSol67] limits the kind of large cardinal embeddings that can exist in a small forcing extension. Here I announce a generalization of this theorem to a broad new class of forcing notions. One consequence is that many of the forcing iterations most commonly found in the large cardinal literature create no new weakly compact cardinals, measurable cardinals, strong cardinals, Woodin cardinals, strongly compact cardinals, supercompact cardinals, almost huge cardinals, huge cardinals, and so on. 

\acknowledgement My research has been supported in part by grants from the PSC-CUNY Research Foundation and from the Japan Society for the Promotion of Science. I would like to thank my gracious hosts here at Kobe University in Japan for their generous hospitality.

Large cardinal set theorists today generally look upon small forcing---that is, forcing with a poset $\P$ of cardinality less than whatever large cardinal $\k$ is under consideration---as benign. This outlook is largely due to the Levy-Solovay theorem \cite[LevSol67], which asserts that small forcing does not affect the measurability of any cardinal. (Specifically, the theorem says that if a forcing notion $\P$ has size less than $\k$, then the ground model $V$ and the forcing extension $V^\P$ agree on the measurability of $\k$ in a strong way: the ground model measures on $\k$ all generate as filters measures in the forcing extension, the corresponding ultrapower embeddings lift uniquely from the ground model to the forcing extension and all the measures and ultrapower embeddings in the forcing extension arise in this way.) Since the Levy-Solovay argument generalizes to the other large cardinals whose existence is witnessed by certain kinds of measures or ultrapowers, such as strongly compact cardinals, supercompact cardinals, almost huge cardinals and so on, one is led to the broad conclusion that small forcing is harmless; one can understand the measures in a small forcing extension by their relation to the measures existing already in the ground model. 

Here in this Communication I would like to announce a generalization of the Levy-Solovay Theorem to a broad new class of forcing notions. 

Historically, the Levy-Solovay theorem addressed G\"odel's hope that large cardinals would settle the Continuum Hypothesis (\CH). G\"odel, encouraged by Scott's \cite[Sco61] theorem showing that the existence of a measurable cardinal implies $V\not=L$, had hoped that large cardinals would settle the \CH\ in the negative. But since one can force the \CH\ to hold or fail quite easily with small forcing, the conclusion is inescapable that large cardinals simply have no bearing whatsoever on the Continuum Hypothesis.

Since that time, set theorists have developed sophisticated tools to combine the two central set theoretic topics of forcing and large cardinals. The usual procedure when forcing with a large cardinal $\k$ whose largeness is witnessed by the existence of a certain kind of elementary embedding $j:V\to M$ is to lift the embedding to the forcing extension $j:V[G]\to M[j(G)]$ and argue that this lifted embedding witnesses that $\k$ retains the desired large cardinal property in $V[G]$. In this way, one is led to consider how a measure $\mu$ in the ground model $V$ can relate to a measure $\nu$ in the forcing extension $V[G]$. The measure $\mu$ may {\df extend} to $\nu$ in the simple sense that $\mu\of\nu$ or it may {\df lift} to $\nu$ when the ground model ultrapower $\jmu:V\to M$ agrees with the larger ultrapower $\jnu:V[G]\to M[\jnu(G)]$ on the common domain $V$. In this terminology, the Levy-Solovay theorem says that after small forcing every measure in the ground model both lifts and extends to a measure in the forcing extension and, conversely, every measure in the extension both lifts and extends a measure in the ground model. 

The truth, however, is that in the large cardinal context most small forcing is, as it were, too small. Rather, one often wants to perform long iterations going up to and often beyond the large cardinal $\k$ in question. With a supercompact cardinal $\k$, for example, one often sees reverse Easton $\k$-iterations along the lines of Silver forcing \cite[Sil71] or the Laver preparation \cite[Lav78]. What we would really like is a generalization of the Levy-Solovay theorem that would allow us to understand and control the sorts of embeddings and measures added by these more powerful and useful forcing notions.

Here, I describe such a generalization. For a vast class of forcing notions, including the iterations I have just mentioned, the fact is that every embedding $j:V[G]\to M[j(G)]$ in the extension that satisfies a mild closure condition lifts an embedding $j:V\to M$ from the ground model. In particular, every measure in $V[G]$ that concentrates on a set in $V$ extends a measure on that set in $V$. From this general fact, I deduce that no forcing of this type can create new weakly compact cardinals, measurable cardinals, strong cardinals, Woodin cardinals, supercompact cardinals, huge cardinals and so on. 

The class of forcing notions for which the theorem applies is quite broad. All that is required is that the forcing admit a {\df gap} at some $\d$ below the cardinal $\k$ in question in the sense that the forcing factors as $\P*\Qdot$ where $\P$ is nontrivial, $\card{\P}<\d$ and $\forces\Qdot$ is $\lted$-strategically closed. (A forcing notion is $\lted$-strategically closed when the second player has a strategy enabling her to survive through all the limits in the game in which the players alternately play conditions to build a descending $(\d+1)$-sequence through the poset, with the second player playing at limit stages.) The Laver preparation, for example, admits a gap between any two stages of forcing. Indeed, in the Laver preparation, the tail forcing is fully directed closed, not merely closed or strategically closed. And the same holds for many of the other reverse Easton iterations one commonly finds in the literature. Moreover, in practice one can often simply preface whatever strategically closed forcing is at hand with some harmless small forcing, such as the forcing to add a single Cohen real, and thereby introduce a gap at $\d=\omega_1$. Further, because $\Qdot$ can be trivial, gap forcing includes all small forcing notions. Examples of useful gap forcing notions are abundant.

An embedding $j:\Vbar\to \Mbar$ is {\df amenable} to $\Vbar$ when $j\restrict A\in\Vbar$ for any $A\in\Vbar$. Let me now state the theorem.

\quiet\theorem Gap Forcing Theorem. Suppose that $V[G]$ is a forcing extension obtained by forcing that admits a gap at some $\d$ below $\k$ and $j:V[G]\to M[j(G)]$ is an embedding with critical point $\k$ for which $M[j(G)]\of V[G]$ and $M[j(G)]^\d\of M[j(G)]$ in $V[G]$. Then $M\of V$; indeed $M=V\intersect M[j(G)]$. If the full embedding $j$ is amenable to $V[G]$, then the restricted embedding $j\restrict V:V\to M$ is amenable to $V$. And if $j$ is definable from parameters (such as a measure or extender) in $V[G]$, then the restricted embedding $j\restrict V$ is definable from the names of those parameters in $V$.

The proof of the theorem is forthcoming in a longer technical paper \cite[Ham$\infty$] (submitted to the Journal of Mathematical Logic and currently available on my web page). Now that I have announced the theorem, what I would like to do in this Communication is illustrate its utility by analyzing the types of measures and embeddings that can exist in a gap forcing extension, ultimately showing in the corollaries below that gap forcing, like small forcing, cannot create various new large cardinals. 

Before proceeding let me record the technical closure fact that from $M=V\intersect M[j(G)]$ it is relatively easy to derive the pair of implications that if $M[j(G)]^\l\of M[j(G)]$ in $V[G]$ then $M^\l\of M$ in $V$ and if $V_\l\of M[j(G)]$ then $V_\l\of M$.

Also, in order to avoid confusion on a subtle point, let me remark that given any embedding $j:V[G]\to\Mbar$ we can let $M=\union\set{j(V_\a)\st\a\in\ORD}$, and it is not difficult to see that $j(G)$ is $M$-generic, that $\Mbar=M[j(G)]$ and moreover that $j\restrict V:V\to M$. Thus, while the statement of the theorem concerns embeddings of the form $j:V[G]\to M[j(G)]$, this form of embedding is fully general.

For those readers who are not completely familiar with the bizarre sorts of embeddings $j:V[G]\to M[j(G)]$ that can exist in a forcing extension, let me stress that in general, quite apart from the question of whether $j$ lifts an embedding from the ground model, one must not presume even that $M\of V$. For example, if $\k$ is a Laver indestructible supercompact cardinal in $V$ and we force to add a Cohen subset $A\of\k$ (by itself, this forcing does not admit a gap below $\k$), then $\k$ remains supercompact in the extension $V[A]$, but any embedding $j:V[A]\to M[j(A)]$ must have $A\in M$ and therefore $M\not\of V$. The point is that the theorem really does identify a serious, useful limitation on the sorts of embeddings that exist in a gap forcing extension. 

In what might be an abuse of terminology, I have stated many of the corollaries below in the slogan form {\it gap forcing creates no new such-and-such kind of large cardinals}. What I mean by this is that if a forcing notion $\P$ admits a gap below a cardinal $\k$ having that large cardinal property in $V^\P$, then $\k$ had the same large cardinal property already in the ground model $V$. That is, I am only making a claim about cardinals above the lowest gap of the forcing. Since many of the most common iterations that admit a gap at all admit a gap that is very low---at $\omega_1$, for example, or near the least inaccessible cardinal---these are often the only cardinals one need consider.

\corollary. Gap forcing creates no new weakly compact cardinals.

\proof The proof of the Gap Forcing Theorem relies on the Key Lemma, appearing in various forms also in \cite[Ham98a], \cite[HamShl98] and \cite[Ham98b], concerning the possibility of certain {\df fresh} sequences in a forcing extension; these are sequences that are not in the ground model though all their initial segments are in the ground model. The Key Lemma asserts that if some forcing has the form $\P*\Qdot$, where $\card{\P}\leq\b$ and $\force\Qdot$ is $\lteb$-strategically closed, then it adds no fresh $\l$-sequences for any $\l$ of cofinality above $\b$. Appealing only to the Key Lemma, let me now argue for this corollary. Suppose that $\k$ is weakly compact in $V[G]$, a forcing extension obtained by forcing with a gap below $\k$. The cardinal $\k$, being inaccessible in $V[G]$, is also inaccessible in $V$. Thus, it remains only to prove that $\k$ has the tree property in $V$. If $T$ is a $\k$-tree in $V$, then by weak compactness it must have a $\k$-branch in $V[G]$. Since every initial segment of this branch is in $V$, it follows by the Key Lemma that the branch itself is in $V$, as desired.\qed

\corollary. Gap forcing creates no new measurable cardinals. 

\proof I will show that every measure on a measurable cardinal in the extension extends a measure in the ground model. It is a standard fact that any ultrapower embedding $j:V[G]\to M[j(G)]$ by a measure $\mu$ in $V[G]$ is closed under $\k$-sequences where $\k=\cp(j)$. If the forcing admits a gap below $\k$, the Gap Forcing Theorem implies that $j\restrict V:V\to M$ is definable from parameters in $V$. Thus, since it is well known that a measure can be recovered from its ultrapower embedding via the equivalence $X\in\mu\iff[\id]_\mu\in j(X)$, it follows that we may recover $\mu\intersect V$ in $V$ from the restricted embedding $j\restrict V$. And since it is easy to see that $\mu\intersect V$ is a measure in $V$, the corollary is proved.\qed

The same argument applies to other kinds of measures, such as supercompactness measures, establishing the following more general fact:

\corollary. After gap forcing, every ultrapower embedding in the extension (by a measure on any set) lifts an embedding from the ground model, and every measure in the extension that concentrates on a set in the ground model extends a measure in the ground model. 

\noindent Of course I am referring here only to embeddings with critical point above the gap. As a caution to the reader, let me stress that the theorem does not say that every ultrapower embedding $j:V[G]\to M[j(G)]$ in the extension is the lift of an ultrapower embedding $j\restrict V:V\to M$ in the ground model. Rather, one only knows that the restricted embedding $j\restrict V$ is definable from parameters in $V$; it is possible to arrange that $j:V[G]\to M[j(G)]$ is the ultrapower by a normal measure on the measurable cardinal $\k$ in the forcing extension while the restricted embedding $j:V\to M$ is some kind of strongness extender embedding, and not an ultrapower embedding at all. 

\corollary. Gap forcing creates no new strong cardinals. 
\ref\Strong

\proof What is true is that if $\k$ is $\l$-strong after forcing with a gap at $\d<\k$, and $\l$ is either a successor ordinal or has cofinality larger than $\d$, then $\k$ was $\l$-strong in the ground model. The reason for making this assumption on $\l$ is that when $j:V[G]\to M[j(G)]$ is the $\l$-strongness embedding induced by a canonical extender for such a $\l$, that is, when $(V[G])_\l\of M[j(G)]$ and $M[j(G)]=\set{j(h)(s)\st s\in (V[G])_\l}$, then it is possible to show that $M[j(G)]^\d\of M[j(G)]$ in $V[G]$, and so the Gap Forcing Theorem applies. Consequently, the restriction $j:V\to M$ is definable from parameters in $V$, and since $V_\l\of M[j(G)]$, the technical closure fact mentioned after the Gap Forcing Theorem implies that $V_\l\of M$. So $\k$ is $\l$-strong in $V$, as desired.\qed

What we actually have is the following:

\corollary. After small forcing $\P$ of size less than $\d$, no further $\lted$-strategically closed forcing $\Q$ can increase the degree of strongness of any cardinal $\k>\d$. 

\proof Suppose that $\k$ is $\l$-strong in $V[g][H]$, the extension by $\P*\Qdot$, and $\k>\d$. In the first case, when $\l$ is either a successor ordinal or a limit ordinal of cofinality above $\d$, the previous corollary shows that $\k$ is $\l$-strong in $V$ and hence also in the small forcing extension $V[g]$. For the second, more difficult case, suppose that $\k$ is $\l$-strong in $V[g][H]$ and $\l$ is a limit ordinal with $\cof(\l)\leq\d$. Let $j:V[g][H]\to M[g][j(H)]$ be a $\l$-strong embedding by a canonical extender, so that $M[g][j(H)]=\set{j(h)(s)\st s\in V[g][H]_\l\and h\in V[g][H]}$. Thus, $j$ is the embedding induced by the extender $$E=\set{\<A,s>\st A\of V_\k\and s\in j(A)\and s\in V[g][H]_\l},$$ which is a subset of $P(V_\k)\cross V[g][H]_\l$. This extender is the union of the smaller extenders $E\restrict\b=E\intersect(P(V_\k)\cross V[g][H]_\b)$ for unboundedly many $\b<\l$. By the result of the previous corollary, we may assume that these smaller extenders each extend a strongness extender in $V$. Since each of these extenders extends uniquely to $V[g]$, the small forcing extension, it follows by the strategic closure of $\Qdot$ that $E\intersect V[g]$ is in $V[g]$ and hence $\k$ is $\l$-strong in $V[g]$, as desired.\qed

The two previous results are complicated somewhat by the intriguing possibility that small forcing could actually increase the degree of strongness of some cardinal. This question, an unresolved instance of the Levy-Solovay theorem, is raised in \cite[HamWdn]. One could ask the corresponding question replacing small forcing with gap forcing, {\it is it possible that forcing with a gap below $\k$ can increase the degree of strongness of $\k$?}\/ But the truth of the matter is that the previous corollary shows that if gap forcing $\P*\Qdot$ can increase the degree of strongness of a cardinal, then this increase is entirely due to the initial small forcing factor $\P$. And the only way this can occur is if a $\ltl$-strong cardinal is made $\l$-strong for some limit ordinal $\l$ of small cofinality.

\corollary. Gap forcing creates no new Woodin cardinals. 

\proof If $\k$ is Woodin in $V[G]$, then by definition this means that for every $A\of\k$ there is a cardinal $\g<\k$ that is $\ltk$-strong for $A$, meaning that for every $\l<\k$ there is an embedding $j:V[G]\to M[j(G)]$ with critical point $\g$ such that $A\intersect\l=j(A)\intersect\l$. Such an embedding can be found that is $(\l+1)$-strong and induced by the canonical extender, so we may assume that $M[j(G)]$ is closed under $\g$-sequences. Thus, for $A$ in the ground model, the Gap Forcing Theorem shows that the restricted embedding $j:V\to M$ witnesses the $\l$-strongness of $\g$ for $A$ in $V$, and so $\k$ was a Woodin cardinal in $V$, as desired.\qed 

Define that a forcing notion is {\df mild} relative to $\k$ when every set of ordinals of size less than $\k$ in the extension has a name of size less than $\k$ in the ground model. For example, the reverse Easton iterations one often finds in the literature are generally mild because the tail forcing is usually sufficiently distributive, and so any set of ordinals of size less than $\k$ is added by some stage before $\k$. Additionally, any $\k$-c.c. forcing is easily seen to be mild. 

\corollary. Mild gap forcing creates no new strongly compact cardinals. 

\proof What I mean is that if $\k$ is strongly compact after forcing that admits a gap below $\k$ and that is mild relative to $\k$, then $\k$ was strongly compact in the ground model. Specifically, I will show that if $\k$ is $\theta$-strongly compact after forcing that is mild relative to $\k$ and admits a gap below $\k$, then it was $\theta$-strongly compact in the ground model; and every strong compactness measure in the extension is isomorphic to one that extends a strong compactness measure from the ground model.

The point is that after mild forcing, every strong compactness measure $\mu$ on $P_\k\theta$ in the extension is isomorphic to a strong compactness measure $\tilde\mu$ that concentrates on $(P_\k\theta)^V$. To see why this is so, first notice that since every $\theta$-strongly compact embedding is actually $\theta^\ltk$-strongly compact, we may assume by replacing $\theta$ with $\theta^\ltk$ if necessary that $\cof(\theta)\geq\k$. Now let $j:V[G]\to M[j(G)]$ be the ultrapower by $\mu$, and let $s=[\id]_\mu$. Thus, $j\image\theta\of s\of j(\theta)$ and $\card{s}<j(\k)$. By mildness $s$ has a name in $M$ of size less than $j(\k)$, and using this name we can construct a set $\tilde s\in M$ such that $j\image\theta\of\tilde s\of j(\theta)$ and $\card{\tilde s}<j(\k)$ in $M$. Furthermore, since $\mu$ is isomorphic to a measure concentrating on $\theta$, there must be some ordinal $\z<j(\theta)$ such that $M[j(G)]=\set{j(h)(\z)\st h\in V[G]}$. I may assume that the largest element of $\tilde s$ has the form $\<\a,\z>$, using a suitable definable pairing function, by simply adding such a point if necessary. Let $\tilde\mu$ be the measure germinated by $\tilde s$ via $j$, so that $X\in\tilde\mu\iff\tilde s\in j(X)$. Since $\tilde s$ is a subset of $j(\theta)$ of size less than $j(\k)$ in $M$, it follows that $\tilde\mu$ is a fine measure on $P_\k\theta$ in $V[G]$ that concentrates on $(P_\k\theta)^V$. I will now show that $\mu$ and $\tilde\mu$ are isomorphic. For this, it suffices by the seed theory of \cite[Ham97] to show that every element of $M[j(G)]$ is in the seed hull $\X=\set{j(h)(\tilde s)\st h\in V[G]}\elesub M[j(G)]$ of $\tilde s$. By the choice of $\tilde s$ we know that $\z\in\X$ and so it is easy to conclude that $j(h)(\z)\in\X$ for any function $h\in V[G]$, as desired. So every strong compactness measure is isomorphic to a strong compactness measure that concentrates on $(P_\k\theta)^V$. 

Now the corollary follows because the restricted embedding $j\restrict V:V\to M$ must be definable from parameters in $V$ by the Gap Forcing Theorem, and using this embedding one can recover $\tilde\mu\intersect V$, which is easily seen to be a fine measure on $P_\k\theta$ in $V$, as desired.\qed

\corollary. Gap forcing creates no new supercompact cardinals, almost huge cardinals, huge cardinals, or $n$-huge cardinals for any $n\in\omega$.

\proof If $j:V[G]\to M[j(G)]$ witnesses one of the mentioned large cardinal properties in $V[G]$, then by the Gap Forcing Theorem the restricted embedding $j:V\to M$ is definable from parameters in $V$ and witnesses the very same large cardinal property in $V$.\qed

Let me close with the following observation.

\theorem Observation. The closure assumption on the embedding in the Gap Forcing Theorem cannot be omitted, because if there are two normal measures on the measurable cardinal $\k$ in $V$ then after merely adding a Cohen real $x$ there is an embedding $j:V[x]\to M[x]$ that does not lift an embedding from the ground model.

\proof Suppose that $\mu_0$ and $\mu_1$ are normal measures on $\k$ in $V$ and $x$ is a $V$-generic Cohen real. By the Levy-Solovay Theorem \cite[LevSol67], these measures extend uniquely to measures $\bar\mu_0$ and $\bar\mu_1$ in $V[x]$, and furthermore the ultrapowers by the measures $\bar\mu_0$ and $\bar\mu_1$ in $V[x]$ are the unique lifts of the corresponding ultrapowers by $\mu_0$ and $\mu_1$ in $V$. Let $j:V[x]\to M[x]$ be the $\w$-iteration determined in $V[x]$ by selecting at the $n^\th$ step either the the image of $\bar\mu_0$ or of $\bar\mu_1$, respectively, depending on the $n^\th$ digit of $x$. If $\<\k_n\st n<\w>$ is the critical sequence of this embedding, then for any $X\of\k$ the standard arguments show that $\k_n\in j(X)$ if and only if $X$ is in the measure whose image is used at the $n^\th$ step of the iteration. Suppose now towards a contradiction that the restricted embedding $j\restrict V$ is amenable to $V$. I will show that from $j\restrict P(\k)^V$ one can iteratively recover the digits of $x$. First, by computing in $V$ the set $\set{X\of\k\st \k\in j(X)}$, we learn which measure was used at the initial step of the iteration and thereby also learn the initial digit of $x$. This information also tells us the value of $\k_1=j_{\mu_{x(0)}}(\k)$. Continuing, we can compute in $V$ the set $\set{X\of\k\st \k_1\in j(X)}$ to know the next measure that was used and thereby learn the next digit of $x$ and the value of $\k_2$, and so on. Thus, from $j\restrict P(\k)^V$ in $V$ we would be able to recursively recover $x$, contradicing the fact that $x$ is not in $V$.\qed 

The argument works equally well with any small forcing; one simply uses a longer iteration. 

\medskip
{\parindent=0pt\tenpoint\tightlineskip\sc 
Mathematics, City University of New York, College of Staten Island\par
2800 Victory Boulevard, Staten Island, NY 10314\par 
\tt hamkins@math.csi.cuny.edu\par
http://www.math.csi.cuny.edu/$\sim$hamkins\par}

\quiet\section Bibliography

\nopagenumbers
\parindent=0pt
\newbox\Article
\newbox\Journal
\newbox\Author
\newbox\Vol
\newbox\No
\newbox\Year
\newbox\Page
\newbox\Book
\newbox\Publisher
\newbox\Pubaddr
\newbox\Key
\newbox\Editor
\newbox\Comment
\newbox\Note
\def\entry#1#2\par{\item{#1\quad}\hskip-1.1em#2\par}
\def\article#1{\setbox\Article=\hbox{\sl #1, }}
\def\journal#1{\setbox\Journal=\hbox{\rm #1 }}
\def\author#1{\setbox\Author=\hbox{\sc #1, }}
\def\vol#1{\setbox\Vol=\hbox{\bf #1 }}
\def\no#1{\setbox\No=\hbox{no. #1 }}
\def\year#1{\setbox\Year=\hbox{\rm({\oldstyle #1}) }}
\def\page#1{\setbox\Page=\hbox{\rm p. #1 }}
\def\book#1{\setbox\Book=\hbox{\it #1, }}
\def\publisher#1{\setbox\Publisher=\hbox{\rm #1, }}
\def\pubaddr#1{\setbox\Pubaddr=\hbox{\rm #1, }}
\def\key#1{\setbox\Key=\hbox{#1}}
\def\editor#1{\setbox\Editor=\hbox{\rm(#1, Ed.) }}
\def\comment#1{\setbox\Comment=\hbox{\rm #1}}
\def\note#1{\setbox\Note=\hbox{\rm #1 }}
\def\ref#1\par{\smallskip{#1
  \entry{\ifhbox\Key\unhbox\Key\else[\ ]\fi}%
  \unhbox\Author\unhbox\Note
  \ifhbox\Book \unhbox\Book\unhbox\Publisher\unhbox\Pubaddr
               \unhbox\Editor\unhbox\Page\unhbox\Year\unhbox\Comment
  \else \unhbox\Article\unhbox\Journal\unhbox\Vol\unhbox\No\unhbox\Editor
        \unhbox\Page\unhbox\Year\unhbox\Comment\fi\par}}

\tenpoint\tightlineskip

\ref
\author{Joel David Hamkins}
\article{Canonical seeds and Prikry trees}
\journal{Journal of Symbolic Logic}
\vol{62}
\no{2}
\page{373-396}
\year{1997}
\key{[Ham97]}

\ref
\author{Joel David Hamkins}
\article{Destruction or preservation as you like it}
\journal{Annals of Pure and Applied Logic}
\year{1998}
\vol{91}
\page{191-229}
\key{[Ham98a]}

\ref
\author{Joel David Hamkins}
\article{Small forcing makes any cardinal superdestructible}
\journal{Journal of Symbolic Logic}
\year{1998}
\vol{63}
\no{1}
\page{51-58}
\key{[Ham98b]}

\ref
\author{Joel David Hamkins}
\article{Gap forcing}
\journal{submitted to the Journal of Mathematical Logic}
\comment{(currently available on the author's web page)}
\key{[Ham$\infty$]}

\ref
\author{Joel David Hamkins and Saharon Shelah}
\article{Superdestructibility: a dual to the Laver preparation}
\journal{Journal of Symbolic Logic}
\year{1998}
\vol{63}
\no{2}
\page{549-554}
\key{[HamShl98]}

\ref
\author{Joel David Hamkins and W. Hugh Woodin}
\article{Small forcing creates neither strong nor Woodin cardinals}
\journal{to appear in the Proceedings of the American Mathematical Society}
\key{[HamWdn]}

\ref
\author{Richard Laver}
\article{Making the supercompactness of $\kappa$ indestructible under 
 $\kappa$-directed closed forcing}
\journal{Israel Journal Math}
\vol{29}
\year{1978}
\page{385-388}
\key{[Lav78]}

\ref
\author{Dana S. Scott}
\article{Measurable cardinals and constructible sets}
\journal{Bulletin of the Polish Academy of Sciences, Mathematics}
\vol{9}
\year{1961}
\page{521-524}
\key{[Sco61]}

\ref
\author{Jack Silver}
\article{The Consistency of the Generalized Continuum Hypothesis with the existence of a Measurable Cardinal}
\journal{Axiomatica Set Theory, Proc. Symp. Pure Math. 13 American Mathematical Society}
\vol{I}
\editor{D. Scott}
\year{1971}
\page{383-390}
\key{[Sil71]}

\ref
\author{Levy Solovay}
\article{Measurable cardinals and the Continuum Hypothesis}
\journal{IJM}
\vol{5}
\year{1967}
\page{234-248}
\key{[LevSol67]}

\bye